\documentclass[11pt,leqno]{article}
\usepackage{epsfig}
\usepackage{graphicx}
\usepackage{wrapfig}
\usepackage{amsmath}
\usepackage{amssymb}
\usepackage{amscd }
\usepackage{amsthm}

\usepackage{tikz}
\usetikzlibrary{arrows}

\newcommand{\be}{\begin{eqnarray}}
\newcommand{\ee}{\end{eqnarray}}
\newcommand{\bi}{\begin{itemize}}
\newcommand{\ei}{\end{itemize}}
\newcommand{\bd}{\begin{definition}}
\newcommand{\ed}{\end{definition}}
\newcommand{\bt}{\begin{theorem}}
\newcommand{\et}{\end{theorem}}
\newcommand{\bc}{\begin{corollary}}
\newcommand{\ec}{\end{corollary}}
\newcommand{\bcn}{\begin{conjecture}}
\newcommand{\ecn}{\end{conjecture}}
\newcommand{\br}{\begin{remark}}
\newcommand{\er}{\end{remark}}
\newcommand{\ce}{\begin{eqnarray*}}
\newcommand{\de}{\end{eqnarray*}}
\newcommand{\bpf}{\begin{proof}}
\newcommand{\epf}{\end{proof}}
\newcommand{\bl}{\begin{lemma}}
\newcommand{\el}{\end{lemma}}
\newtheorem{theorem}{Theorem}[section]

\newtheorem{conjecture}[theorem]{Conjecture}
\newtheorem{lemma}[theorem]{Lemma}
\newtheorem{remark}[theorem]{Remark}
\newtheorem{definition}[theorem]{Definition}
\newtheorem{proposition}[theorem]{Proposition}
\newtheorem{Examples}[theorem]{Examples}
\newtheorem{corollary}[theorem]{Corollary}
\numberwithin{equation}{section}

\def\e{\varepsilon}

\def\a{\alpha}

\def\b{\beta}
\def\d{\delta}

\def\l{\lambda}

\def\[{{\Big[}}
\def\]{{\Big]}}
\def\<{{\langle}}
\def\>{{\rangle}}

\def\({{\Big(}}
\def\){{\Big)}}

\def\bt{\begin{theorem}}
\def\et{\end{theorem}}
\def\bl{\begin{lemma}}
\def\el{\end{lemma}}
\def\br{\begin{remark}}
\def\er{\end{remark}}
\def\bx{\begin{Examples}}
\def\ex{\end{Examples}}
\def\bd{\begin{definition}}
\def\ed{\end{definition}}
\def\bp{\begin{proposition}}
\def\ep{\end{proposition}}
\def\bc{\begin{corollary}}
\def\ec{\end{corollary}}
\def\cA{{\mathcal A}}
\def\cB{{\mathcal B}}
\def\cC{{\mathcal C}}

\def\cF{{\mathcal F}}
\def\cG{{\mathcal G}}
\def\cH{{\mathcal H}}

\def\cJ{{\mathcal J}}

\def\cP{{\mathcal P}}

\def\cR{{\mathcal R}}

\def\cW{{\mathcal W}}

\def\cZ{{\mathcal Z}}

\def\mN{{\mathbb N}}

\def\mX{{\mathbb X}}

\def\mZ{{\mathbb Z}}

\usepackage{graphicx}

\begin{document}

\section* {\center{Totally Rank One Interval Exchange Transformations}}
\begin{center}
\textsc
{
Yue Wu,\footnote{Schlumberger WesternGeco Geosolution, Houston, Texas, USA.}
Dongmei Li,\footnote{Corresponding Author; Harbin University of Science and Technology, School of Applied Science, Harbin, Heilongjiang, China.}
Diquan Li,\footnote{Central South University, School of Geosciences and Info-physics,University, Changsha, Hunan, China.}
Yunjian Wang,\footnote{Schlumgerger SWT Technology Center, Houton, Texas, USA.}
}
\end{center}

\begin{abstract}

 For irreducible interval exchange transformations, we study the relation between the powers of induced map and the  induced maps of powers and raise a condition of equivalence between them. And skew production of Rauzy induction map is set up and verified to be ergodic regarding to a product measure. Then we prove that almost all the interval exchange transformations are totally rank one (rank one for all powers of positive integers) by interval. As a corollary, for almost all interval exchange transformations, rank one transformations are dense $G_{\d}$ in the weak closure.\\ \\
 {\bfseries Keywords}  \, Rank One, Interval Exchange Transformations, Rauzy-Veech Induction\\
 {\bfseries 2000 MR Subject Classification}  \, 37A30

\end{abstract}
\section* {}

In this paper, we extend Veech's theorem about rank one interval exchange transformations (i.e.t.) (\cite{KEA}, \cite{VIA}) to all positive-integer powers of interval exchange transformations (totally rank one). Property of rank one is a crucial conception in the field of ergodic theory and dynamic systems. Originated from stacking construction of transformations related to symbolic dynamics, rank one transformation has been studied with regarding to many topics (\cite{FER1}, \cite{GLA}, \cite{KIN1}), including functional spectrum, invariant measures of productions, and many other topics. If the base of each Rohlin tower of the rank one transformation is an interval, we say it is rank one by interval.If all the powers of a transformation are rank one, we say it is totally rank one (definition \ref{Total}). Here it is showed that measure theoretically all the interval exchange transformations are totally rank one by interval (Theorem \ref{Main}).

\bd [Totally Rank One] \label{Total}
We say a finite measure-preserving system $(\mX, \b, T, \mu)$ is totally rank one, if for every $q\in \mN$, $T^{q}$ is rank one.
\ed

In Section one, we introduce the fundamental concepts of i.e.t, Rauzy-Veech induction (\cite{RAU}, \cite{VEE1}, \cite{VEE3}) and Veech's theorems (\cite{VEE3}) about i.e.t.'s related to this paper. At the end of section one, the main theorem of this paper is stated: almost all interval exchange transformations are totally rank one by interval. In Section 2, a condition for the induced map of a power of i.e.t. to be equal to the power of a induced map of i.e.t. is given, with redar to Rauzy-Veech induction. In Section 3, we introduce a skew product, which is extended from the Rauzy map. It is showed that this skewing transformation is ergodic and conservative relative to a product measure. In the last section, Section 4, the major propositions studied in Section 2 and Section 3 are used to set up a cutting and stacking structure, related to the Rauzy-Veech induction, of a power of i.e.t. This shows that measure theoretically, all the i.e.t.'s are totally rank one by interval. A corollary is conducted that rank one transformations are dense $G_{\d}$ in the weak closure of a i.e.t, measure theoretically.

\section{Introduction}\label{section1}

Let $\Lambda _{m}\subset R^{m}$ be a positive cone, and $\lambda=(\lambda_{1},\cdots ,\lambda _{m})\in \Lambda_{m}$.
Let $\cG_{m}$ be the group of m-permutations, and
$\mathcal{G}_{m}^{0}$ be the subset of $\mathcal{G}_{m}$ which
contains all the irreducible permutations on $\{1,2,\cdots, m\}$. A
permutation $\pi$ is irreducible if and only if for any $1\leq k
<m,\,\{1,2,\cdots,k\}\neq\{\pi (1),\cdots,\pi(k)\}$, or
equivalently $\sum\limits_{j=1}^{k}(\pi(j)-j)>0,\, (1\leq k<m))$. Given $\lambda \in \Lambda_{m},\, \pi \in \mathcal{G}_{m}^{0}$,
the corresponding interval exchange transformation is defined by:\\
 \begin{equation}\begin{array}{l}
 T_{\l,\pi}(x)=x-\beta_{i-1}(x)+\beta_{\pi i-1}(\lambda^{\pi}),\quad (x\in
[\beta_{i-1}(\l),\quad
\beta_{i}(\l))\,),\\

\begin{array}{l}
\text{where }  \\
              \quad \lambda^{\pi}=(\l_{\pi^{-1}1},\l_{\pi^{-1}2},\cdots,
\l_{\pi^{-1}m}).\\
               \quad \beta_{i}(\lambda)=\left\{\begin{array}{clcr} 0& i=0\\
\sum\limits_{j=1}^{i}\l_{j}& 1\leq i\leq m .
                      \end{array}\right.
 \end{array}
 \end{array}
\end{equation}
Obviously $\beta_{\pi i-1}(\l^{\pi})=\sum\limits_{j=1}^{\pi
i-1}\l_{\pi^{-1}j}$,
 and the transformation $T_{\l,\pi}$, which is also denoted by $(\l,\pi)$, sends the
$i$th interval to the $\pi (i)$th position.\\

{\em\bf Rauzy-Veech induction.\/\/} For $T_{\l,\pi}$, the Rauzy map
sends it to its induced map on
$[0,\left|\l\right|-min\left\{\l_{m},\l_{\pi^{-1}m}\right\})$,
which is the largest admissible interval of form
$J=[0,L),0<L<\left|\l\right|$.\\
Given any permutation, two actions $a$ and $b$ are:
\be
a(\pi)(i)=\left\{\begin{array}{llcr}\pi(i) & i\leq\pi^{-1}m
\\\pi(i-1)&\pi^{-1}m+1<i\leq m \\ \pi(m)
&i=\pi^{-1}m+1 \end{array}\right.
\ee
 and
 \be
b(\pi)(i)=\left\{\begin{array}{llcr}\pi(i) & \pi(i)\leq\pi(m)
\\\pi(i)+1&\pi(m)+1<\pi(i)< m \\ \pi(m)+1
&\pi(i)=m . \end{array}\right.
\ee
 \indent The Rauzy-Veech map $\mathcal{Z}(\l,\pi):\,
\Lambda_{m}\times \mathcal{G}_{m}^{0}\rightarrow \Lambda_{m}
\times \mathcal{G}_{m}^{0}$ is determined by :
\be\label{RzVch}\mathcal{Z}(\l,\pi)=(A(\pi ,c)^{-1}\l,c\pi) ,
\ee
where $c=c(\l,\pi)$ is defined by
\be
c(\l,\pi)=\left\{\begin{array}{clcr}a,&\l_{m}<\l_{\pi^{-1}m}\\
b,&\l_{m}>\l_{\pi^{-1}m.}\end{array}\right.
\ee
$\mathcal{Z}(\l,\pi)$ is a.e. defined on
$\Lambda_{m}\times\left\{\pi\right\}$, for each $\pi\in\mathcal{G}_{m}^{0}$. \\
\indent The matrices $A=A(\pi,c)$ in \ref{RzVch} are defined as the following:
\be \label{MatrixA} A(\pi ,a)=\left(\begin{tabular}{c|c}
$I_{\pi^{-1}m}$&$\begin{array}{ccccc}
0&0&\cdots&0&0\\
0&0&\cdots&0&0\\
.&.&\cdots&.&.\\
0&0&\cdots&0&0\\
1&0&\cdots&0&0
  \end{array}$\\\hline\\
0& $\begin{array}{ccccc}
0&1&\cdots&0&1\\
0&0&\cdots&0&0\\
.&.&\cdots&.&.\\
0&0&\cdots&0&1\\
1&0&\cdots&1&0
  \end{array}$\\
\end{tabular}\right)
\ee

\be \label{MatrixB} A(\pi ,b)=\left(
\begin{tabular}{c|c}
 $I_{m-1}$ &0\\\hline\\
$\underbrace{\begin{array}{ccccccc}
0&\cdots&0&1&0&\cdots&0\end{array}}_{\mbox{1 at the jth position}}$&1\\
\end{tabular}
\right)
\ee
where $I_{k}$ is the $k$-identity matrix, and $j=\pi^{-1}m$.\\
\indent And the normalized Rauzy map $\cR : \, \Delta_{m-1}\times \cG^{0}_{m}\rightarrow
 \Delta_{m-1}\times\cG^{0}_{m}$  is defined by
 \be
 \cR (\l, \pi)=(\frac{\displaystyle A(\pi ,c)^{-1}\l}{\displaystyle \left|A(\pi ,c)^{-1}\l\right|},
 c\pi)=(\frac{\displaystyle \pi^{*}_{1}\cZ (\l
,\pi)}{\displaystyle \left|\pi^{*}_{1}\cZ (\l ,\pi)\right|},
\pi^{*}_{2}\cZ (\l ,\pi)),
 \ee
 where $\pi^{*}_{1}$ and $\pi^{*}_{2}$ are the projection to the first coordinate and the second coordinate respectively.\\
Iteratively,
\be\label{Cn}
\mathcal{Z}^{n}(\l,\pi)=((A^{(n)})^{-1}\l,
c^{(n)}\pi)=(\l ^{(n)},\pi^{(n)}),
\ee
where
\be\label{An}
c^{(n)}=c_{n}c_{n-1}\cdots c_{1},(c_{1},\cdots,c_{n}\in
 \{a,b\}, c_{i}=c(\mathcal{Z} ^{i-1}(\l,\pi)))
\ee
and
\be
A^{(n)}=A(\pi,c_{1})A(c^{(1)}\pi, c_{2})A(c^{(2)}\pi,c_{3})\cdots
 A(c^{(n-1)}\pi,c_{n}) .
 \ee
\indent The Rauzy class $\cC\subseteq\cG_{m}$ of $\pi$ is a set of orbits for the group of maps generated by $a$ and $b$. On the $\cR$ invariant component $\Delta_{m-1}\times\cC$, we have:
\bt [H.Masur\cite{MASUR};W.A. Veech\cite{VEE1}]\label{Ergodic} Let $\pi\in\cG_{m}^{0}$, the set of irreducible permutations. For Lebesgue almost all $\l\in\Lambda_{m}$, normalized Lebesgue measure on $I^{\l}$ is the unique invariant Borel probability measure for $T_{(\l,\pi)}$. In particular, $T_{(\l,\pi)}$ is ergodic for almost all $\l$.
\et
\bt [W.A. Veech\cite{VEE3}]\label{Rkone} Let $\pi\in\cG_{m}^{0}$, the set of irreducible permutations. For Lebesgue almost all $\l\in\Lambda_{m}$, normalized Lebesgue measure on $I^{\l}$ is rigid and rank one, thus admits simple spectrum.
\et\
We extend Veech's result to the following theorem:
\bt \label{Main} All the notations as above,
suppose $m>1$, $m\in \mN$, $\pi\in \cG^{*}_{m}$, $q\in \mN$. For
Lebesgue almost all $\lambda\in \Lambda$, $T^{q}_{(\lambda, \pi)}$
is rank one by interval with flat stack, thus it is also rigid.
\et

The theorem will be seen to imply that for the corresponding
corresponding interval exchange transformations, each has a residual
set of  rank one transformations in its weak closure, by [KIN], also
in its commutant. \\

Next we introduce some result which
would be utilized in section 5
and section 6:
\bt \label{ErgInt} Let $m>1$, $\pi\in
\cG^{0}_{m}$, $\cR =\cR (\pi)$ is the Rauzy class containing
$\pi$. Then there exists on $\triangle_{m-1}\times\cR (\pi)$ a
unique absolutely continuous invariant measure $\mu$ up to a
scalar mutiple, such that $\cP$ (see (2.1.9) ) is conservative and
ergodic on $\triangle_{m-1} \times \cR$ relative to $\mu$. The
density of $\mu$ on $\triangle _{m-1} \times \{\pi\}$ is the
restriction of a function which is positive, rational, and
homogeneous of degree
$m$ on $\Lambda_{m}\times \{\pi\}$.
\et
Next,the induced map of $\cP$ on $\triangle_{m-1}\times
\{\pi\}$ will be described, and a condition for a measurable
function to be
essential constant will be given. \\

By Theorem \ref{ErgInt}, induced map $\cP _{\pi}$ of $\cP$ on $\triangle_{m-1}\times
\{\pi\}$ is well defined measure theoretically. Then for a.e. $x=(\lambda ,\pi)\in
\triangle_{m-1}\times\{\pi\}$ there exists $n(x)\in\mN$ such that
$\cP _{\pi}(x)=\cP^{n(x)}(x)$, and there exists $A=A(x)$ such that
$\cP_{\pi}(x)=(\frac{\displaystyle A^{-1}\lambda}{\displaystyle
\left|A^{-1}\lambda\right|}, \pi)$. What's more, if we let
$\triangle^{A}=(A\Lambda_{m})\cap \triangle_{m-1}$, and define
$P_{A}:\triangle ^{A}\rightarrow
\triangle_{m-1}$ by

\be
P_{A}(\alpha)=\frac{\displaystyle A^{-1}\alpha}{\displaystyle
\left|A^{-1}\alpha\right|}\quad (\alpha\in \triangle^{A})
\ee

then for any $y=(\alpha ,\pi)\in \triangle^{A}\times \{\pi\}$, $\cP_{\pi}(y)=(P_{A}(\alpha),\pi)$.\\
There exists a set $\cF$ containing countable many elements (which
are visitation matrices),such that $\cF_{1}=\{\triangle ^{A}|A\in
\cF\}$ is a partition of $\triangle _{m-1}$ into infinitely many atoms.\\
 \\
$\; \mathbf {Note.}$:  For any $A\in \cF$, for Lebesgue almost all
$\lambda\in \triangle ^{A}, \cP^{i}_{\pi}(\lambda,\pi)$ is defined for all $i\in \mN$.\\
\\

To end this section, we recall (9.5) of [VEE3]
\bt \label{eigenConst}
Let $\cH$ be a separable
Hilbert space, and $U_{A}$ be unitary operators on $\cH$. Suppose
$f:\triangle _{m-1}\rightarrow \cH$ is a measurable function and
$f(\cP_{\pi}x)=U_{A}f(x)$, for $a.e. \;x\in \triangle^{A}, A\in
\cF)$, then $f$ is constant.
\et

\section{Powers  of the Induced Map and the Induced Map of
the powers}\label{section1}

 Without normalization, the $n$th iteration of the Rauzy map raises an m-interval exchange
 trasformation on a subinterval $J_{n}$ with vector
 $\alpha=\lambda^{(n)}$  and m-permutation $\pi^{(n)}$ (see \ref{Cn}). Note that
 $\left|J_{n}\right|=\left|\alpha\right|$. In this section we pay
 more attention on the first return time of $T_{\lambda,
 \pi}|_{J_{n}}$. The visitation matrix $A^{(n)}$ shows the
 orbit distribution of $I^{(\alpha)}_{j}$ among $I^{(\lambda)}_{ i}$ 's. That
 is $A_{ij}^{(n)}$ is the number of the visitation of $T^{l}(I_{\alpha j})$ ($0 \leq l <
 a_{j}, a_{j}$ is the first return time of $I^{(\alpha)}_{ j}$ to $J_{n}$) on $I_{i}$.
 Thus the first return time of any point in $I^{(\alpha)}_{ j}$ is
 $a_{j}=\sum\limits_{i=1}^{m}A_{ij}^{(n)}$, that is $a_{j}$ equals the
 $j$th column sum of $A^{(n)}$. Canonically, this yields a
 stack structure with m-stacks, $S_{1},S_{2},\cdots ,S_{m}$, each corresponding to the $T$
 orbits of a subinterval of $\alpha$.
 The $j$-th stack $S_{j}$ has base $I^{(\alpha)}_{ j}$ and height
 $a_{j}$. $T$ sends each level (except for the top) of $S_{j}$ to
 the higher level, and sends the top back into $J_{n}=\cup I^{(\alpha)}_{ j}$.
 Suppose the first return time of $T^{q}$ to $J_{n}$ is $r_{q}(x), x\in J_{n}$,
 then $T^{q r_{q}(x)}(x)\in J_{n}$, $T^{q l}(x)\notin J_{n}(l=1,2,\cdots,
 r_{q}(x)-1)$. It is easily seen that though $T^{qr_{q}(x)}\in J_{n}$, the
 return time of $T$ to $J_{n}$ may be strictly less than
 $qr_{q}(x)$.\\

\bl \label{EvenRetn}
$\quad$All notations as above,
suppose $a^{(k)}_{j}$ is the $j$th column sum of $A^{(k)}$, and
each $a_{j}^{k}$  is odd, then there exists $1\leq j_{0}\leq m$
such that $a_{j_{0}}^{(k+1)}$ is even.
\el
\begin{proof}
Suppose $\cZ^{k}(\lambda, \pi)=(\lambda^{(k)},\pi^{(k)})$, let
$\pi_{0}=\pi^{(k)}$, then $A^{(k+1)}=A^{(k)}A$, where $A$ is either \ref{MatrixA}
or \ref{MatrixB}\\
\\
For the first case we see that \\
$$a^{(k+1)}_{\pi ^{-1}_{0} m+1}=a^{(k)}_{\pi^{-1}_{0}m}+a^{(k)}_{m}$$\\
for the second case we see that \\
$$a^{(k+1)}_{\pi^{-1}_{0}m}=a^{(k)}_{\pi^{-1}_{0}m}+a^{(k)}_{m}$$
\end{proof}
\bd
$T_{(\lambda,\pi)}$
satisfies the infinite distinct orbit condition (i.d.o.c., by Kean[\cite{KEA}]) if  ${T^{-j}(\b_{i})}_{j\in \mN}$, $1\leq i \leq m$ (the negative trajactories of discontinuities of $T$) are infinite disjoint sets¡£
\ed
\bl
$\quad$ If $T_{(\lambda,\pi)}$
satisfies the i.d.o.c, then $T^{q}_{(\lambda, \pi)}$ also
satisfies the i.d.o.c.
\el
\begin{proof}
Suppose $T^{q}_{(\lambda, \pi)}$ is $T_{(\xi,\tau)}$, then the
$q(m-1)$ discontinuities are $D_{q}=\{T^{-i}\beta_{i},1\leq j< m,
0\leq i < q\}$. One may concern the point 0. Since
$0=T(\beta_{k})$ for some $1\leq k\leq m$, the $i$ th preimage of
$0$ is the $i-1$th preimage of $\beta_{k}$, already contained in
$D_{q}$. So what need to be shown is that the $T^{q}$ orbit of
$T^{-i_{1}}\beta_{j}$ and that of $T^{-i_{2}}\beta_{j}$ are
disjoint. This is true since $i_{1}\neq i_{2}$ and $0\leq i_{1},
i_{2}<q$.
\end{proof}
\bp \label{InEqRetn1}
$\quad$If there exist $j_{1},
j_{2}\in \mN$ , $1\leq j_{1},j_{2}\leq m$ such that $a_{j_{1}}\neq
a_{j_{2}} \mod q$, then $((T^{q}|_{J_{n}})\circ(T|_{J_{n}})\neq
(T|_{J_{n}})(T^{q}|_{J_{n}}))$.
\ep
\begin{proof}
Let $\Gamma=\{i|a_{i}\equiv a_{j_{1}}\mod q\}$
\\
then $\Gamma^{*}=\{i|a_{i}\neq a_{j_{1}} \mod q\}$ since
$a_{j_{1}}\neq a_{j_{2}}\mod q$, we know that $\Gamma\neq \phi$,
$\Gamma^{*}\neq \phi$, and $j_{1}\in \Gamma ,\quad j_{2}\in
\Gamma^{*}$. Correspondingly, let
\\
$$K=\{x|x\in I^{(\alpha)}_{ i}, i\in\Gamma \}$$
$$K^{*}=\{x|x\in I^{(\a)}_{ i}, i\in\Gamma^{*} \}$$

 Then since $T^{q}_{\lambda,\pi}$ satisfies i.d.o.c.,
$T^{q}|_{J_{n}}(K)\cap K^{*}\neq \phi$, there exits $x\in
T^{q}|_{J_{n}}(K)\cap K^{*}$. Since $T^{q}|_{J_{n}}$ is
invertible, there exits $y\in K\cap (T^{q}|_{J_{n}})^{-1}(K^{*})$.
Thus

$$(T^{q}|_{J_{n}})\circ (T|_{J_{n}}(y))=T^{ql_{1}+a_{j_{1}}}(y),\quad l_{1}\in
\mZ$$
$$(T|_{J_{n}})\circ (T^{q}|_{J_{n}}(y))=T^{ql_{2}+a_{j_{0}}}(y),\quad l_{2}\in \mZ,j_{0}\in \Gamma
^{*}$$

Since $a_{j_{1}}\neq a_{j_{0}}\quad \mod q$\\
$$(T^{q}|_{J_{n}})\circ (T|_{J_{n}})(y)\neq (T|_{J_{n}})\circ
(T^{q}|_{J_{n}})(y)$$
\end{proof}

\begin{remark}
Proposition \ref{InEqRetn1} shows that if two
subintervals of the induced map have different return times modulo
$q$, then the induced map of $T$ and that of $T^{q}$ on the same
interval $J_{n}$ do not commute.
\end{remark}
From Lemma \ref{EvenRetn} and Proposition \ref{InEqRetn1}, we can see  that
during the consequent iterations of Rauzy map, at least one of
$J_{2k}, J_{2k+1}$ is such an subinterval such that
$T^{2}|_{J_{k'}}$ and $T|_{J_{k'}}$ do not commute, $k'\in \{2k,
2k+1\}$. One may draw more general conclusion about the general case
of $T^{q}\quad (q>1, q\in \mN)$. From these propositions, relation
between $T^{q}|_{J_{n}}$ and $T|_{J_{n}}$ are studied. If we equalize the degree and study
$T^{q}|_{J_{n}}$ and $(T|_{J_{n}})^{q}$, it is interesting that
these two transformation with the same domain and range (both are
$J_{n}$ ) may be equal to each other. A condition for it to be true
is given right after this paragraph, one will see it is
crucial for proving Theorem \ref{Main}\\

Suppose $a_{i}\equiv 1 \mod q$, $i=1,2,\cdots , m$; then for any
$x\in J_{n}$ there exists a sequence of positive integers $q$:
$a_{i_{1}}(x),a_{i_{2}}(x),\cdots ,a_{i_{q}}(x)$ such that
$T^{a^{*}(x)}=(T|_{J_{n}})^{q}(x)$,
$a^{*}(x)=\sum\limits^{k=1}_{q}a_{i_{k}}(x)$ and
$T^{a_{i_{j}}(x)}((T^{q}|_{J_{n}})^{j-1}(x))=(T^{q}|_{J_{n}})^{j}(x)$,
$1\leq j\leq q $ . That is to say that $\sum\limits^{j}_{k=1}
a_{i_{k}}(x)$ is the $j$th return time of $x\in J_{n}$ under the
transformation $T$. Since $a_{i}\equiv 1 \mod q$,
$i=1,2,\cdots,m$, $a^{*}\equiv 0 \mod q$. We claim that \\
\be \label{commute1}
T^{q}|_{J_{n}}=(T|_{J_{n}})^{q}
\ee

Suppose $(T^{q}|_{J_{n}})(x)=T^{a'(x)}$, to prove \ref{commute1},  the only
thing needs to be verified is that $a'(x)\geq a^{*}(x)$. ($a'(x)\leq
a^{*}(x)$ is impled by $a^{*}\equiv 0 \mod q)$ . On the other hand,
the $j$th return time of $x\in J_{n}$ under $T$ is
$a^{*}_{j}(x)=\sum\limits^{j}_{k=1}a_{ik}(x)$, $a^{*}_{j}(x)\equiv
 j \neq 0 \mod q $, $(1\leq j<q)$, thus $a'(x)\geq a^{*}(x)$. That
 proves the following theorem:
\bt
All notations and
definitions as above, if $a_{i}\equiv 1 \mod q$,
$T^{q}|_{J_{n}}=(T|_{J_{n}})^{q}$. \begin{center} {\bf $\S 2$ An
Ergodic Skew Product of Rauzy map \/}
\end{center}
\et

As an extension of the Rauzy map, a skewing transformation will be
defined. In this section, we will show that this skew product is
ergodic and conservative relative to a certain product measure.
This result is the fundament tool to prove the major result,
Theorem \ref{Main}. We locate the case to that
the skewing group is a finite group, and we reach the result a little bit indirectly, i.e.
we prove the ergodic property of a skew product conjugate to that we want.\\

Starting with $m,q\in \mN ,m,q\geq 2$, a finite group $G$ is given
by $G=GL(m,\mZ _{q})$, where $\mZ _{q}$ is the ring of integers
modulo $q$. A map $g:\cG^{0}_{m}\times \{ a,b \} \rightarrow G $
is defined by:
\be
g(\pi, c)=A(\pi,c) \mod q,\,(\pi \in \cG^{0}_{m},c\in \{a,b\})
\ee

The same notation $g$ is used for a map $g: X\rightarrow G$, where
$X=\triangle_{m-1}\times \cR(\pi)$, $\cR(\pi)$ the Rauzy class of
$\pi$. Since they are naturally associated, it is well understood.
That is $g(x)=g(\pi(x),c(x))$, $c(x)=c(\lambda (x),\pi (x))$.\\
\\
The normalized Rauzy map is $\cP (\lambda ,
\pi)=(\frac{\displaystyle A^{-1}\lambda }{\displaystyle \left|
A^{-1} \lambda\right|},c(\pi))$, $\cP _{\pi}$ is the induced map
of $\cP$ on $\triangle _{m-1} \times \{\pi\}$. Suppose for $n\in
\mN$, $\cP ^{n}_{\pi}(x)=\cP ^{r}(x)$, where $r=r_{1}+r_{2}+\cdots
+r_{n}$, and if $\overline{r_{i}}=r_{1}+r_{2}+\cdots +r_{i}$, then
$\cP ^{i}_{\pi}(x)=\cP ^{\overline{r_{i}}}(x)$, $g^{(r)}(x)=g(\cP
^{r-1}(x))g(\cP ^{r-2}(x))\cdots g(\cP (x))g(x)$ is associated
with $x$. One can understand $g^{(r)}(x)$ as a closed path from
$\pi $ to $\pi$. Let $G'(\pi)=\left\{ g^{(r)}(x) | x\in X(\pi),
\cP^{r}x\in X(\pi) \right\}$. Next, it is shown that $G'(\pi)$ is
a semigroup of $G$ based on the fact that $\cP _{\pi}$ is an
expansion.
\bt
$G'(\pi)$ is closed
under multiplication.
\et
\begin{proof}
Suppose $g^{(r)}(x),g^{(s)}(y)\in G'(\pi),\quad x,y\in X(\pi)$.
Suppose $\cP^{s}(y)$ is associated with the visitation matrix $A$
then $y\in\triangle^{A}\times \{\pi\}$. We know that
$\cP^{s}(\triangle^{A}\times\{\pi\})=\triangle_{m-1}\times
\{\pi\}$. Therefore there exists $y_{0}\in X(\pi)$ such that
$g^{(s)}(y)=g^{(s)}(y_{0})$ and $\cP^{s}(y_{0})=x$. Then
$$g^{(r)}(x)g^{(s)}(y)=g^{(r)}(\cP^{s}(y_{0}))g^{(s)}(y_{0})$$
$$\qquad=g^{(r+s-1)}(y_{0})\cdots g^{(s)}(y_{0})g^{(s-1)}(y_{0})\cdots g(y_{0})$$
$$\qquad=g^{(r+s)}(y_{0})$$

Thus $$g^{(r)}(x)g^{(s)}(y)\in G'(\pi)$$\\
\end{proof}

$G'(\pi)$ is thus a sub-semigroup of $G$. let $G(\pi)=G'(\pi)$ ,
then $G(\pi)$ is a subsemigroup of the finite group $G$. Therefore $G(\pi)$ is a subgroup of $G$.\\

If we have two permutations $\pi_{1}, \pi_{2}$ in the same
irreducible Rauzy class $\cR (\pi)$, by Theorem \ref{Ergodic} (Veech's
ergodic theorem), there exists $x\in X(\pi_{1})$ and $n\geq 0$
such that $\cP^{n}(x_{1}) \in X(\pi_{2})$. Assign $g^{(n)}(x)$ to
$g(\pi_{1},\pi_{2})$, of course there may be other choice, but
once a fixed one is assigned, it is well defined.
$g(\pi_{1},\pi_{2})$ understood as a path from $\pi_{1}$ to
$\pi_{2}$, two 'paths' may be connected in the following sense: if
$\pi_{1},\pi_{2},\pi_{3}\in \cR(\pi)$, $y\in X(\pi_{2})$,
$\cP^{s}(y)\in X(\pi_{3})$ and $g^{s}(y)=g(\pi_{2},\pi_{3})$, then
there exists $x\in X(\pi_{1})$ such that $\cP^{r}(x)=y$ and
$g^{r}(x)=g(\pi_{1},\pi_{2})$, then
$g^{r+s}(x)=g^{s}(y)g^{r}(x)g(\pi_{2},\pi_{3})g(\pi_{1},\pi_{2})$.
It is obvious that $g(\pi_{2},\pi_{1})g(\pi_{1},\pi_{2})\in G'(\pi_{1})$.\\
Now we are ready to define a skewing map $\cW _{0}:X\times
G_{0}\rightarrow X\times G_{0}$ where $G_{0}=G(\pi_{0})$ for some
fixed $\pi_{0}\in \cR$. That is:
$$\cW_{0}(x,\gamma)=(\cP x, h(x)\gamma)$$
where $$h(x)=g^{-1}(\pi_{0}, c(x)\pi (x))g(x)g(\pi_{0},\pi(x))$$

  To understand $\cW_{0}$, we do some computation
$$\cW_{0}^{2}=(\cP^{2}x, h(\cP x)h(x) \gamma)$$
$$\cdots$$
$$\cW_{0}^{r}=(\cP ^{r}_{x},h(\cP^{r-1} x)h(\cP^{r-2} x)\cdots h(\cP x)h(x))$$

Let $$h^{(r)}(x)=h(\cP^{r-1} x)h(\cP^{r-2} x)\cdots h(\cP x)h(x)$$

Next define a space $X^{*}$ and a transformation $\cW$ on $X^{*}$:
$$X^{*}=\underset{\pi\in \cR (\pi_{0})}{\cup} X(\pi)\times g(\pi_{0},\pi)G_{0}$$
$$\cW(x,g(\pi_{0},\pi (x))\gamma)=(\cP x, g(x)g(\pi_{0},\pi (x) )\gamma)$$

The relation between $\cW$ and $\cW_{0}$ is given by $\varphi
:X^{*}\rightarrow X\times G_{0}$ defined as
$\varphi(x,g(\pi_{0},\pi (x))\gamma)=(x,\gamma)$, it is easy to
see that $\varphi$ is bijective. Since it is also true that
$\varphi \cW=\cW_{0} \varphi$, we have $\cW$, $\cW_{0}$ are conjugate.
\bl \label{compose1}
If $\pi\in \cR$, $c\in\left\{
a, b\right\}$, then $G(c\pi)=g(\pi , c)G(\pi)g(\pi, c)^{-1}$.
\el
\begin{proof}
Since $$g(\pi , c)G(\pi)g(\pi,c)^{-1}\subseteq G(c\pi)$$
$$g(\pi , c)^{-1}G(c\pi)g(\pi,c)\subseteq G(\pi)$$

We have $$g(\pi , c)g(\pi,
c)^{-1}G(c\pi)g(\pi,c)g(\pi,c)^{-1}\subseteq g(\pi, c)G(\pi)g(\pi,
c)^{-1}\subseteq G(c\pi)$$

Thus $$G(c\pi)=g(\pi , c)G(\pi)g(\pi ,
c)^{-1}$$
\end{proof}

\bl \label{Struct1}
Let $H(\pi)=\left\{h^{(r)}(x)|x\in X(\pi), \cP ^{r}(x)\in
X(\pi)\right\}$, Then $H(\pi)=G_{0}$.
\el
\begin{proof}
It is obvious that $g^{r}(x)$ generates $G(\pi)$ , since $G(\pi)$
is finite, $\left\{g^{r}(x)|x\in X(\pi), \cP^{r}(x)\in X(\pi)\right\}=G'(\pi)$.\\

By Lemma \ref{compose1} we have $G(\pi)g(\pi_{0},\pi)\gamma=g(\pi_{0},\pi)G_{0}$. Since $G(c\pi)$ and $G(\pi)$ are conjugate, $G(\pi_{1})$ and $G(\pi_{2})$ are
conjugate for any $\pi_{1},\pi_{2}\in \cR(\pi_{0})$. Thus $G(\pi)$
and $G_{0}$ have the same number of elements.
Since $G (\pi)g(\pi_{0},\pi)\gamma \subset g(\pi_{0},\pi) G_{0}$, $G(\pi)g(\pi_{0},\pi)\gamma=g(\pi_{0},\pi)G_{0}$.\\

Because $\varphi$ is a conjugation,
$\varphi(\cW^{r}(x,g(\pi_{0},\pi)\gamma))=\cW^{r}_{0}\varphi(x,g(\pi_{0},\pi)\gamma)$,
it is obvious that $\# H(\pi)=\# (G(\pi))=\# (G_{0})$. Since
$H(\pi)\subset G_{0}, H(\pi)=G_{0}$.

\end{proof}
Next we show the ergodic property of the two skewing
transfromations.
\bl \label{MeasureSkew}
All notations as above,
and suppose $\omega \otimes \delta_{\pi}$ is the invariant measure
of $\cP_{\pi}$, $\sigma$ is the normalized Haar measure on
$G_{0}$. Then $\cW_{0}$ is ergodic and conservative relative to
the measure $\nu_{0}=\sum\limits_{\pi \in \cR} \omega \otimes
\delta_{\pi}\otimes \sigma$.
\el
\begin{proof}
First the parallel properties of the induced map on $X(\pi)\times
G_{0}$ is verified, then those of $\cW_{0}$ itself will be gained.\\

Denote the induced map by $\cW_{\pi}$, that is
$\cW_{\pi}(x,\gamma)=(\cP_{\pi}x, h_{\pi}(x)\gamma)$, where
$\cP_{\pi}(x)=\cP^{r}(x)$ and $h_{\pi}(x)=h^{(r)}(x)$.\\

Suppose $F\in L^{2}(\omega \otimes \delta_{\pi} \otimes \sigma)$
is essentially invariant under $\cW_{\pi}$, that is $F\circ
\cW_{\pi}(z)=F(z)$ a.e. $z$. Then define a function
$f:X(\pi)\rightarrow \cH=L^{2}(G_{0})$ by
$f(x)(\gamma)=F(x,\gamma)$.\\

$\;$Recall the partition $\cF_{1}$ from Section 0. Suppose
$A=A^{(n)}(\lambda, \pi), (\lambda \in \triangle)$, let $\triangle
=\triangle ^{A} \in \cF_{1}$ . It is also true that $\triangle$ is
associated with common $c_{1},c_{2}, \cdots ,c_{n}$ and
$\pi_{1},\pi_{2}, \cdots ,\pi_{n}$,
where $c_{i}=c(\cP ^{i-1} x)$, $\pi_{i}=\pi(\cP ^{i-1} x)$.\\

Since $$h^{(n)}(x)=h(\cP^{n-1}x) \cdots h(x)$$
$$=g^{-1}(\pi_{0}, c(\cP^{n-1}x)\pi (\cP^{n-1} x))g(\cP^{n-1} x)g(\pi_{0}, \pi (\cP^{n-1}(x))$$
$$\cdots g^{-1}(\pi_{0}, c(x)\pi(x)) g(x) g(\pi_{0},\pi (x))$$
$$=g^{-1}(\pi_{0}, c_{n}(x)\pi_{n}(x))g(x)g(\pi_{0},\pi_{n}x) $$
$$\cdots g^{-1} (\pi_{0}, c(\pi_{1})) g(x) g(\pi_{0}, \pi_{1})$$
$h^{n}(x)$ is constant on $\triangle \times \left\{\pi \right\}$.
Thus a unitary operator $U_{\triangle}: \cH \rightarrow \cH$ may
be defined as
$U_{\triangle}(\psi)(x)=\psi((h^{(n)}(x))^{-1}\gamma)$, $\psi \in
\cH , x\in \triangle \times \left\{\pi \right\}$. Since $\sigma$
 is the Haar measure(counting measure here) on $G_{0}$, $U_{\Delta}$ is an isometry, thus a unitary operator.\\

Since $F\circ \cW _{\pi}(z)=F(z)$, $ z\in X_{\pi}\times
G_{0}$, we have:
\begin{align}\label{eigen1}
\nonumber
f(\cP_{\pi}x)(\gamma) & =F(\cP_{\pi}x,\gamma)=F(\cW_{\pi}^{-1}(\cP_{\pi}x,\gamma))\\
& =F(x,(h_{\pi})^{-1}\gamma )=U_{\triangle}(f(x))
\end{align}

By Theorem \ref{eigenConst},
we know that $f$ is essentially constant. This tells us: \\

$\qquad$ {\em i)\/} Let
$F(\gamma)=F(x_{1},\gamma)=F(x_{1},\gamma)$ for such a $x\in
X(\pi)$, it is well understood that $F(\gamma)$ is a fixed
function in $\cH=L^{2}(\sigma)$.\\

$\qquad$ {\em ii)\/} From \ref{eigen1} we know that $F(h^{(n)}(x)\gamma)=F(\gamma)$\\

Now we apply Lemma \ref{Struct1} , $H(\pi)\gamma =G_{0}$, thus $F(\gamma)$
is constant on $G_{0}$. That is to say if $F_{1}$ is any $\cW_{0}$
invariant function in $L^{2}(X\times G_{0})$, the restricted map
of $F_{1}$ on $X(\pi)\times G_{0}$ is constant. Since $\cW_{0}$
send an element in $X(\pi)\times G_{0}$ to an element in
$X(c(\pi))\times G_{0}$ and $\cR (\pi)$ is the set
$\left\{c_{1}c_{2}\cdots c_{k}(\pi), k\geq 0, c_{i}\in \left\{ a,
b\right\}, 1\leq i\leq k\right \}$, it is not hard to see that
$F_{1}$ is constant on $X\times G_{0}$. Therefore $\cW_{0}$ is
ergodic. $\cW_{0}$ is also conservative, since $\cP_{\pi}$ is
conservative for all $\pi \in \cR (\pi)$, $\cR (\pi)$ and $G_{0}$
are both finite.
\end{proof}
\bt \label{ErgodicSkew}
$\cW$ is ergodic and
conservative relative to the measure $\nu = \sum \limits _{\pi \in
\cR}\omega \otimes \delta_{\pi} \otimes g(\pi_{0}, \pi) \sigma$.
\et
\begin {proof}
Since $\cW$ and $\cW_{0}$ are conjugate, Lemma \ref{MeasureSkew} implies
this result.
\end{proof}

\begin{center} {\bf $\S 3$ Totally Rank One Property for Interval Exchange Transformations\/}
\end{center}

In this section, $\pi \in G^{*}_{m}$ is fixed. That is to say an
interval exchange is associated with $(\lambda , \pi)$ ($\lambda \in
\triangle _{m-1}$). Based on theorem \ref{ErgodicSkew}, we shall show the
measure theoretic generic rank one property of
the powers of the $\pi$-interval exchange transformations(Theorem \ref{Main}).\\

A positive matrix can be associated with the iteration of the Rauzy
transformation on
$$X_{\pi}=\left\{(\lambda, \pi)| \lambda \in
\Delta_{m-1}, \pi \mbox{ a given m-permutation}\right\}$$
One way to
see this is that once $(\lambda, \pi)$ satisfies i.d.o.c, it is
minimal, the iterations are corresponding to the induced maps of
$(\l ,\pi)$ on smaller and smaller subintervals. What's more,
Theorem \ref{ErgodicSkew} implies that, for a.e. $\lambda$, the orbit of $(\lambda,
\pi)$ will visit $U\times\{\pi\}\times \{ e\}$ infinitely many
times, for any open subset of $U\subset \Delta_{m-1}$. If $M$ is any $m\times m$ positive matrix, then $M\cdot
A(\pi^{*}, c)$ $(\pi^{*}\in \cR(\pi))$ is also a positive matrix.
Therefore, there exists $c_{1}, c_{2}, \cdots, c_{n}$ with each
$c_{i}\in \{a, b\}$ $1\leq i\leq n$ such that $c_{n}c_{n-1}\cdots
c_{1} \pi=\pi$, $B=A^{(n)}$ is positive and
$B\equiv e\mod q$. $A^{(n)}$ may be expressed by induction:\\
I
$$A^{(1)}=A(\pi, c_{1})$$
II \\

{\em 3.3.1\/}
\be
A^{(i+1)}=A^{(i)}A(c_{i}\cdots
c_{2}c_{1}(\pi), c_{i+1})
\ee
$\;$It is nice to see that for any $\eta \in B \Lambda$, $\cZ
^{n}(\eta,\pi)=(B^{-1} \eta ,\pi )$ passing the same sequence of
permutations as
$c_{1}\pi, c_{2}c_{1} \pi , \cdots , c_{n}c_{n-1} \cdots c_{1}\pi $. \\

Given $\tau >0$, define an open set in $\triangle_{m-1} \times \{
\pi \}\times \{e\}$ ($e$ be the identity in $G_{0}$), that is $\cJ
=\{(\alpha , \pi , e)\}| \alpha \in \triangle _{m-1},
\alpha_{1}>1-\tau$. Let $\cJ ^{*}=\{(\frac{\displaystyle B\alpha}{\displaystyle \left|B \alpha\right|}, \pi, e)| (\alpha , \pi , e)\in \cJ\}$. Then $\omega ^{*}(\cJ^{*})>0$.\\
The above process will be continued in the proof of the following lemma:
\bl \label{Fund}
$\quad$All notations as above, for a.e.
$\lambda\in \Delta_{m-1}$ and $\e >0$ there is an positive integer
$n_{0}$ and an interval $J\in [0, \left| \lambda\right|)$ such that:\\
(a) $T^{iq} J$ are pairwise disjoint, $0\leq i <n_{0}$.\\
(b) $T^{q}$ is linear on $T^{iq}J$, $0\leq i<n_{0}$\\
(c) $\sum\limits^{n_{0}-1}_{i=0}\left|T^{iq} J\right|>(1-\epsilon)\left|\lambda\right|$\\
(d) $\left|J\cap T^{n_{0}q}J\right|>(1-\epsilon)\left|\lambda\right|$
\el

\begin{proof}
 We prove the corresponding results a.e. $\lambda \in \triangle_{m-1}$, equivalently.\\

Assume $\lambda \in \triangle_{m-1}$ satisfies: there exists
$k\in\mN$ such that $\cW^{k}(\lambda, \pi, e)\in \cJ^{*}$, based
on Theorem \ref{ErgodicSkew}. In other words such $\l$ form a set full measure.\\

Since $\cW^{k}(\l ,\pi ,e)$ is an element in $\cJ^{*}$,
$\cW^{k+n}(\l ,\pi ,e)$ is an element in $\cJ$. Suppose
$\cW^{k}(\l ,\pi , e)$ is associated with the visitation matrix
$A$. Then $A\equiv e \mod q$. Suppose $\cW^{k+n}(\lambda,
\pi,e)=(\overset{\sim}{\alpha},\pi,e)$ then the associated matrix
is $A^{(k+n)}(\lambda, \pi)=AB$, $AB$ is positive and $AB\equiv e
\mod q$. Also we know that
$\overset{\sim}{\alpha}=\frac{\displaystyle
(AB)^{-1}\lambda}{\displaystyle \left|(AB)^{-1}\lambda\right|}$.
Let $\alpha=(AB)^{-1}\lambda$, $J'=[0, \left|\alpha\right|)$, then
since $\overset{\sim}{\alpha}_{1}>(1-\tau)$, it is true that
$\alpha_{1}>(1-\tau)\left|\alpha\right|$. \\
\\
 $\;$Suppose $(T|_{J'})^{q}=T_{(\eta , \varsigma)}$ with $\eta
\in \Lambda_{q(m-1)+1}$, $\left|\eta\right|=\left|\alpha \right|$,
 $\varsigma \in \cG^{0}_{q(m-1)+1}$.\\
 \\

By induction, it is easy to see that there exists $1\leq k\leq
q(m-1)+1$ such that $\eta_{k}>(1-2^{q-1}\tau)\left|\alpha\right|$.
Let $J$ be the $k$-th interval of $\eta$. Since $AB\equiv e \mod
q$, we have $(T|_{J'})^{q}=(T^{q}|_{J'})$ by Theorem 3.1.6. Now
visiting each $I^{(\l)}_{j}$ the same number of times, points in
each suninterval of $\eta$ come back to $J'$ coincidently for the
first $q$ times under $T$. Therefore, if $I^{\eta}_{i}$ is the
$i$-th interval of $\eta$, we can let $a^{(t)}_{i}$ be the first
return time to $J'$ of $(T|_{J'})^{t-1}(I^{\eta}_{i})$ under $T$.
Let $\overset {\sim}{a_{i}}=\sum\limits^{q}_{t=1}
a^{(t)}_{i}=a^{*}_{i}q$.Thus
$$T^{q}|_{J'}(I_{\eta_{i}})=(T^{q})^{a^{*}_{i}}(I_{\eta_{i}})$$
Let $n_{0}=a^{*}_{k}$, then \\

$\qquad$$\cA)$ $\quad \cdots\quad(T^{q})^{l}(J)$, $(0\leq l< n_{0})$ are pairwise disjoint and\\

$\qquad$$\cB)$ $\quad \cdots\quad\left|(T^{q})^{n_{0}}(J)\cap J\right|>(1- 2^{q}\tau)\left| \alpha \right|$. \\

Given an positive $m\times m$ matrix $M$, define
$$v(M)=\max\limits_{1\leq i,j,k\leq m}\frac{M_{ij}}{M_{ik}}$$
then we have (by [VEE2]) if $m_{j}=\sum\limits^{m}_{i=1} M_{ij}$,
then
$$m_{j}\leq v(M)m_{k}\quad (1\leq j, \, k\leq m)$$
$$v(PM)\leq v(M),\quad \mbox{P a nonnegative} \quad m \times m \quad \mbox{matrix}$$

Therefore, $$v(AB)\leq v(B)$$ \\

  Recall that each $a^{(t)}_{i}$ is one of the column summations of
$AB$. Thus $$a^{(t)}_{i}\leq v(B) a^{(s)}_{j},1\leq i, j\leq
q(m-1)+1,\quad1\leq t,s \leq q$$
$\;$This implies $a^{*}_{i}\leq v(B)a^{*}_{j}$, $1\leq i, j \leq q(m-1) +1$.\\
\\
$\;$So the total length(measure) of the remaint of the major stack
($\cup^{n_{0}-1}_{l=0}(T^{q})^{l}(J)$) is:
$$\left|\lambda \right| -n_{0}\left|J \right|= \left| \sum\limits^{m(q-1)+1}_{i=1} a^{*}_{i} \eta_{i}\right|-a^{*}_{k} \eta_{k}$$
$$=\sum\limits ^{m(q-1)+1}_{i=1, i\neq k}a^{*}_{i}\eta _{i}\leq v(B)a^{*}_{k}2^{q-1}\tau \left| \eta\right|$$
$$\leq v(B)a^{*}_{k}\frac{2^{q-1}\tau}{1-2^{q-1}\tau}\eta_{k}\leq v(B)\frac{2^{q-1}\tau}{1-2^{q-1}\tau} $$
That is \\

  $\qquad$$\cC)$ $\quad \cdots\quad\left|\cup^{n_{0}-1}_{l=0}(T^{q})^{l}(J)\right| >(1-v(B)\frac{2^{q-1}\tau}{1-2^{q-1}\tau})$, $(\left|\lambda\right|=1)$.\\
In order to make the approximation coincide to $\e >0$ choose $\tau$
small enough such that $2^{q}\tau <\epsilon$ and
$v(B)\frac{2^{q-1}\tau}{1-2^{q-1}\tau}<\epsilon$. Combining
$\cA)\quad \cB )\quad \cC )$ and Theorem \ref{ErgodicSkew} (for a fixed
$\pi\in\cG^{0}_{m}$ almost all $\l$, $(\l ,\pi)$ will visit
$\cJ^{*}$ infinitely often under $\cW$), we abtain (a)(b)(c)(d) in
the Lemma.
\end {proof}

It is easy to see that Theorem \ref{Main} is a corollary of Lemma \ref{Fund}.
\bp \label{Base}
Let $(\mX , \cB , \mu, T)$ be an
automorphism system on a standard measure space. Suppose all powers
of $T$, such that this subset is a dense $G_{\delta}$ set in
$Wcl(T)$, and all transformations in it are rank one.
\ep
\begin{proof}
 We say two partitions $P_{1}$ and $P_{2}$ satisfying
$P_{1}>_{\varepsilon} P_{2}$ if for any atom $p\in P_{2}$,
there exists $p_{1}, p_{2}, \cdots , p_{m} \in P_{1}$ such that $m(p\Delta (\cup ^{m}_{i=1}p_{i}))<\varepsilon$.\\

Next we use $P^{(n)}$ to denote the partition by dyadic sets of rank $n$.\\

Let $R(n, q, \varepsilon)=\{S|$ there exists $B \subset [0, 1)$,
such that $P_{B,q}>_{\varepsilon} P^{n}$, where $P_{B,q}=\{B, SB ,\cdots , S^{q-1}B , X-\cup^{q-1}_{i=0}S^{i}B\}\}$\\

Then $R(n,q,\epsilon)$ is an open set in $G=Aut (\mX ,\cB ,\mu)$.\\

Suppose $T$ has all powers $(T^{n}, n\in \mN)$ rank one and rigid. Then $Wcl(T)$ has uncountable many elements.
$W(n, q, \varepsilon)= Wcl(T)\cap R(n,q,\varepsilon)$ is an open set
in $Wcl (T)$.\\

Since all powers of $T$ are rank one, we have
$$T^{n}\in \cup _{q}W(n, q, \varepsilon),n\in \mN$$
and $\{T^{n}, n\in \mN\}$ are dense in $Wcl(T)$.Thus $\cup _{q}W(n, q, \varepsilon)$ is a dense open set in $Wcl(T)$.
Therefore $\cap_{n}\cup _{q}W(n, q, \varepsilon)$ is a dense $G_{\delta}$ set
in $Wcl(T)$, with all elements rank one.
\end{proof}

\bc
   All the notations as above,
suppose $m>1$, $m\in \mN$, $\pi\in \cG^{*}_{m}$, $q\in \mN$. For
Lebesgue almost all $\lambda\in \Lambda$, it is true that there is a
dense $G_{\delta}$ subset of rank one transformations in the weak
closure of $T^{q}_{(\lambda, \pi)}$.
\ec
\begin{proof}
By theorem \ref{Main} and proposition \ref{Base}.
\end {proof}

\end{document}